\documentclass[a4paper,12pt]{article}
\RequirePackage{amsfonts,amssymb}
\begin{document}
\author{P.Bieliavsky$^*$, M.Pevzner\footnote{Research supported by 
the Communaut\'e fran\c caise de Belgique, through an
Action de Recherche Concert\'ee de la Direction de la Recherche 
Scientifique and the grant NWO 047-008-009}\\{\small \tt pbiel@ulb.ac.be, mpevzner@ulb.ac.be}
\\ Universit\'e Libre de
Bruxelles
\\Belgium}
\title{Symmetric Spaces and Star representations II : Causal Symmetric Spaces }
\date{  }
\maketitle
\addtolength{\topmargin}{-1.5cm}
\newtheorem{thm}{Theorem}[section]
\newtheorem{lem}[thm]{Lemma}
\newtheorem{englth}[thm]{Theorem}
\newtheorem{pro}[thm]{Proposition}
\newtheorem{cor}[thm]{Corollaire}
\newtheorem{de}[thm]{Definition}
\newtheorem{rmq}[thm]{Remarque}

          \newtheorem{dfn}{Definition}[section] 
          \newtheorem{rmk}{Remark}[section]
          \newtheorem{prop}{Proposition}[section] 
          \newtheorem{exs}{Examples}[section] 
          \newcommand{\D}{\mbox{$\cal D$}}
          \newcommand{\BR}{\mbox{$\Bbb R$}} 
          \newcommand{\BC}{\mbox{$\Bbb C$}}
          \newcommand{\Levi}{\mbox{$\cal L$}} 
          \newcommand{\Lg}{\mbox{$\frak g$}}
          \newcommand{\Lk}{\mbox{$\frak k$}} 
          \newcommand{\La}{\mbox{$\frak a$}} 
          \newcommand{\Lh}{\mbox{$\frak h$}}
          \newcommand{\Lm}{\mbox{$\frak m$}} 
          \newcommand{\Ln}{\mbox{$\frak n$}}
          \newcommand{\Lt}{\mbox{$\frak t$}} 
          \newcommand{\Ll}{\mbox{$\frak l$}}
          \newcommand{\Lz}{\mbox{$\frak z$}} 
          \newcommand{\LH}{\mbox{$\frak H$}}
          \newcommand{\Lu}{\mbox{$\frak u$}} 
          \newcommand{\Lb}{\mbox{$\frak b$}}
          \newcommand{\Lc}{\mbox{$\frak c$}} 
          \newcommand{\Ls}{\mbox{$\frak s$}}
          \newcommand{\Lr}{\mbox{$\frak r$}} 
          \newcommand{\Lq}{\mbox{$\frak q$}}
          \newcommand{\Le}{\mbox{$\frak e$}} 
          \newcommand{\Lhr}{\mbox{${\frak h}_r$}}
          \newcommand{\Lqr}{\mbox{${\frak p}_r$}} 
          \newcommand{\Lhs}{\mbox{${\frak h}_s$}}
          \newcommand{\Lqs}{\mbox{${\frak p}_s$}}
          \newcommand{\Lsu}{\mbox{$\frak{su}$}}
          \newcommand{\Lso}{\mbox{$\frak{so}$}} 
          \newcommand{\Lsp}{\mbox{$\frak{sp}$}}
          \newcommand{\Lspin}{\mbox{$\frak{spin}$}}
          \newcommand{\Lgl}{\mbox{$\frak{gl}$}}
          \newcommand{\ad}{\mbox{ad}}
          \newcommand{\adk}{\mbox{$\mbox{ad}_{\Lk}$}}
          \newcommand{\adh}{\mbox{$\mbox{ad}_{\Lh}$}}
          \newcommand{\adx}{\mbox{$\mbox{ad}_{\xi}$}} 
          \newcommand{\pr}{\mbox{pr}}
          \newcommand{\Pf}{{\em Proof}. }
          \newcommand{\EPf}{\hfill$\Box$}
          \newcommand{\holm}{\frak{hol}} 
          \newcommand{\rad}{\mbox{rad}}
          \newcommand{\bxi}{\mbox{$\bar{X}_i$}} 
          \newcommand{\bxj}{\mbox{$\bar{X}_j$}}
          \newcommand{\bxk}{\mbox{$\bar{X}_k$}} 
          \newcommand{\byi}{\mbox{$\bar{Y}_i$}}
          \newcommand{\byj}{\mbox{$\bar{Y}_j$}} 
          \newcommand{\byk}{\mbox{$\bar{Y}_k$}}
          \newcommand{\tOmega}{\mbox{$\tilde{\Omega}$}}
          \newcommand{\tnabla}{\mbox{$\tilde{\nabla}$}}
          \newcommand{\hOmega}{\mbox{$\hat{\Omega}$}}
          \newcommand{\OO}{{\mbox{${\cal O}$}}}

%          \renewcommand{\theenumi}{\(\alph{enumi}\)} 
%\setcounter{section}{-1}
%%% d/but
\let\Dbl\Bbb
\let\BBox\Box
\def\Box{$\BBox$}
%\font\goth=eufm10
\def\g{\mathfrak{g}}
\def\k{\mathfrak{k}}
\def\t{\mathfrak{t}}
\def\s{\mathfrak{s}}
\def\u{\mathfrak{u}}
\def\z{\mathfrak{z}}
\def\a{\mathfrak{a}}
\def\n{\mathfrak{n}}
\def\p{\mathfrak{p}}
\def\P{\mathfrak{P}}
\def\h{\mathfrak{h}}
\def\gl{\mathfrak{gl}}
\def\e{\mathfrak{e}}
\def\l{\mathfrak{l}}
\def\V{\mathfrak{V}}
\def\W{\mathfrak{W}}

\def\C{\mathbb{C}}
\def\R{\mathbb{R}}
\def\N{\mathbb{N}}
\def\Z{\mathbb{Z}}
\def\qed{\hfill\quad\raise -2pt \hbox{\vrule\vbox to 10pt
{\hrule width 8pt
\vfill\hrule}\vrule}\newline}

\def\Det{\mathop{\rm Det}\nolimits}
\def\et{\mathop{\rm et}\nolimits}
\def\id{\mathop{\rm id}\nolimits}
\def\Re{\mathop{\rm Re}\nolimits}
\def\Exp{\mathop{\rm Exp}\nolimits}

\def\tr{\mathop{\rm tr}\nolimits}

\def\mod{\mathop{\rm mod}\nolimits}

\def\Ind{\mathop{\rm Ind}\nolimits}

\def\Tr{\mathop{\rm Tr}\nolimits}

\def\rank{\mathop{\rm rank}\nolimits}

\def\diag{\mathop{\rm diag}\nolimits}

\def\Int{\mathop{\rm int}\nolimits}
\def\carre{\mathbin{\hbox{\vrule\vbox to 4pt
    {\hrule width 4pt \vfill\hrule}\vrule}}}

\def\Aut{\mathop{\rm Aut}\nolimits}

\def\ch{\mathop{\rm ch}\nolimits}

\def\tanhyp{\mathop{\rm th}\nolimits}

\def\Exp{\mathop{\rm Exp}\nolimits}

\def\ad{\mathop{\rm ad}\nolimits}

\def\Im{\mathop{\rm Im}\nolimits}

\def\Ad{\mathop{\rm Ad}\nolimits}

\def\dim{\mathop{\rm dim}\nolimits}

\def\Sym{\mathop{\rm Sym}\nolimits}

\def\spur{\mathop{\rm spur}\nolimits}

\def\div{\mathop{\rm div}\nolimits}

\ifx\optionkeymacros\undefined\else\endinput\fi

\abstract{\footnotesize We construct and identify
star representations canonically associated with holonomy reducible simple
symplectic symmetric spaces. This leads the a non-commutative 
geometric realization of the correspondence between causal symmetric 
spaces of Cayley type and Hermitian symmetric spaces of tube type.\newline MSC2000: 22E46, 53C35, 81S10}
\section*{ Introduction}
The aim of this paper is twofold. We first want to present a
generalization
of a construction given in \cite{[Biel]}. There, a covariant star product (see 
Definition~\ref{star-cov}) has been defined on a dense open subset of every 
holonomy reducible simple symplectic (non-Kaehler) symmetric space 
$M=G/H$, providing a star representation $\rho$ of the transvection 
algebra $\g$ of the symmetric space $M$ (see Section 2). Despite the 
fact that this star product is only defined on an open subset of $M$, 
the representation $\rho$ in the case $\g=\s\l(2,\R)$ exponentiates to the
group $G=SL(2,\R)$ as a holomorphic discrete series representation. In 
particular the covariant star product does {\sl not} lead to a 
representation of $SL(2,\R)$ prescribed by the orbit type---indeed, in
this 
case, one has $M=SL(2,\R)/SO(1,1)$ which classically yields a principal 
series representation. It should also be noted that the star 
representation does generally not exponentiate to $G$ as in the case of 
$G=SL(3,\R)$ for instance. In the present work, we prove that the star 
representation $\rho$ exponentiates to $G$ when the symmetric space 
$M=G/H$ is of Cayley type. The resulting representation of $G$ turns out
to 
be a holomorphic discrete series representation (when one assigns a 
particular real value to the deformation parameter); see Theorem \ref{resultat}. 
The proof uses Jordan techniques, providing close relations between Jordan 
algebras theory and covariant star products.\\
Second, we would like to indicate how the above construction leads to a 
(non-commutative) realization of the (algebraic) duality existing between 
Cayley symmetric spaces and Hermitian symmetric spaces of tube type. 
Below, we make this assertion more precise.\\
Let $G/K$ be a Hermitian symmetric space of tube type. Then there exists 
one and only one (up to isomorphism) symmetric space $G/H$ such that
\begin{enumerate}
\item $H$ acts reducibly on  $T_{eH}(G/H)$,
\item $G/H$ carries a  $G$-invariant symplectic structure.
\end{enumerate}
This fact leads to the well-known ``duality" between Hermitian 
symmetric spaces of tube type and Cayley symmetric spaces \cite{[Faraut Olafsson]}. 
As 
such, this duality is algebraic in the sense that it is a correspondence 
between two lists of involutive simple Lie algebras. Another duality 
defined on Hermitian symmetric spaces is the so called
``compact-non-compact" 
duality. In this case, the correspondence is not only algebraic. Indeed, 
denoting by $U$ a compact real form of $G^\C$, there is a holomorphic 
$G$-equivariant embedding $G/K\to U/K$ underlying the duality. Equivalently, 
one has a homomorphism of algebras $C^\infty(U/K)\to C^\infty(G/K)$.
Back to Cayley symmetric spaces, our construction of a covariant star 
product described above leads to a deformation of the (infinitesimal) action 
of (the 
Lie algebra) of $G$ on a open subset of $G/H$. This deformed 
action turns out to be
equivalent to the action of $G$ on the tube domain $G/K$. Analogously 
with the compact-non-compact duality, one therefore gets a geometric 
realization of the correspondence between Cayley symmetric spaces and 
Hermitian symmetric spaces of tube type. 

This paper is organized as follows. In Sections 1 and 2, we recall the 
notion of covariant star product \cite{[Arnal]} and results of
\cite{[Biel]}. Section 3 contains our main result, it starts with 
recalling some Jordan algebras theory.

\section{Covariant $\star$-products}
In this section, $(M,\omega)$ is a symplectic manifold and $\g$ is a 
finite dimensional real Lie algebra. One assumes one has a 
representation of $\g$ as an algebra of symplectic vector fields on 
$(M,\omega)$. That is, one has a Lie algebra homomorphism 
$\g\rightarrow{\cal X}(M):\:X\mapsto X^{*}$ (${\cal X}(M)$ stands for 
the space of smooth sections of $T(M)$) such that for all $X$ in $\g$
one has
$$
{\cal L}_{X^{*}}\omega=0,
$$
where $\cal L$ denotes the Lie derivative. One supposes furthermore 
that this representation of $\g$ is strongly Hamiltonian which means 
that there exists a $\Bbb R$-linear map 
$$
\g\stackrel{\lambda}{\mapsto}C^{\infty}(M):\:\:X\mapsto\lambda_{X},
$$
such that
\begin{eqnarray*}
(i)&&d\lambda_{X}=i_{X^{*}}\omega,\:\forall X\in\g;\\
(ii)&&\lambda_{[X,Y]}=\{\lambda_{X},\lambda_{Y}\},
\end{eqnarray*}
where $\{\:,\}$ denotes the Poisson structure on $C^{\infty}(M)$ 
associated to the symplectic form $\omega$.

\begin{de}\label{strham} A quadruple $(M,\omega,\g,\lambda)$ with $(M,\omega),\g$ 
and $\lambda$ as above is called a strongly hamiltonian system. The 
map
$\lambda:\g\mapsto C^{\infty}(M)$ is called the moment mapping.
\end{de}
{\bf Example 1.} Coadjoint orbits. Let $M={\cal O}\subset\g^{*}$ be a 
coadjoint orbit of a Lie group $G$ with Lie algebra $\g$. In this case,
we denote by $\g\rightarrow{\cal X}(M):\:X\mapsto X^{*}$ the rule
which associates to an element $X$ in $\g$ its {\it fundamental vector
field} on $\cal O$ :
$$
X_{x}^{*}:=\frac d{dt}\vert_{0}\Ad^{*}(\exp(-tX))x,
$$
where $\Ad^{*}(g)x$ denotes the coadjoint action of the element $g\in 
G$ on $x\in\g^{*}$.

The formula $\omega_{x}^{\cal O}(X^{*},Y^{*}):=\langle 
x,[X,Y]\rangle$ (with $X,Y\in\g$) then defines a symplectic 
structure called after, Kirillov, Kostant and Sauriau, the {\it KKS 
symplectic
form} on $\cal O$. Defining, for all $X\in\g$, the function 
$\lambda_{X}\in C^{\infty}(\cal O)$ by
$$
\lambda_{X}(x):=\langle x, X\rangle,
$$
one then has that the quadruple $({\cal O},\omega^{\cal 
O},\g,\lambda)$
is a strongly Hamiltonian system.

In the setting of deformation quantization there is a natural way to 
define the quantization of a classical hamiltonian system \cite{[Arnal]},\cite{[BFFLS]},\cite{[F]}. We first recall what
deformation quantization  (star product) is.
\begin{de}
Let $(M,\{\:,\:\})$ be a Poisson manifold. A star product on 
$(M,\{\:,\:\})$
is an associative multiplication $\star_{\nu}$ on the space 
$C^{\infty}(M)[[\nu]]$ of formal power series in 
the parameter $\nu$ 
with coefficients in the smooth complex-valued functions on $M$. One 
furthermore requires the following properties to be true.

(i) The map $\star_{\nu}:C^{\infty}(M)[[\nu]]\times 
C^{\infty}(M)[[\nu]]\rightarrow 
C^{\infty}(M)[[\nu]]$ is $\Bbb C[[\nu]]$-bilinear and for all $u\in 
C^{\infty}(M)\subset C^{\infty}(M)[[\nu]]$ one has 
$u\star_{\nu}1=1\star_{\nu}u=u$; ($\Bbb C[[\nu]]$ denotes the field 
of power 
series in $\nu$ with (constant) complex coefficients).

(ii) For all $u,v\in C^{\infty}(M)$ one has,
\begin{eqnarray*}
1.&& u\star_{\nu}v \mod(\nu)=uv,\\
2.&& (u\star_{\nu}v-v\star_{\nu}u)\mod (\nu^{2})=2\nu\{u,v\}.
\end{eqnarray*}
\end{de}
In other words, a star product is an associative formal deformation 
of 
the pointwise multiplication of functions in the direction of the 
Poisson structure.

{\bf Example 2.} The Moyal star product on $\Bbb R^{2n}$. We fix 
$M=\Bbb R^{2n}$ (or an open set in $\Bbb R^{2n}$) and $\omega 
=\sum_{i<j}\Lambda_{ij}dx^i\wedge dx^j$ with constant coefficients 
$\Lambda_{ij}$'s.

Let $u,v\in C^{\infty}(\Bbb R^{2n})$ and define their Moyal product 
by the formal power series:
$$
u\star_{\nu}^{0}v:= 
uv+\sum_{k=1}^{\infty}\frac{\nu^k}{k!}\Lambda^{i_{1}j_{1}}\ldots
\Lambda^{i_{k}j_{k}}\frac{\partial^k}{\partial 
x^{i_{1}}\ldots\partial 
x^{i_{k}}}u\:\frac{\partial^k}{\partial x^{j_{1}}\ldots\partial 
x^{j_{k}}}v.
$$
The $\Bbb C[[\nu]]$-bilinear extension of the above product to 
$C^{\infty}(\Bbb R^{2n})[[\nu]]$ then defines a star product on 
$(\Bbb R^{2n},\Lambda)$ called Moyal star product.
\begin{de}\label{star-cov}
\cite{[Arnal]} Let $(M,\omega,\g,\lambda)$ be a classical strongly 
hamiltonian system. A star product $\star_{\nu}$ on $(M,\omega)$ is 
called
$\g$-covariant if for all $X,Y\in\g$ one has
$$
\lambda_{X}\star_{\nu}\lambda_{Y}-\lambda_{Y}\star_{\nu}\lambda_{X}=2\nu\{\lambda_{X},
\lambda_{Y}\}.
$$
\end{de}
In order to avoid technical difficulties in defining the star 
representation (see below), we will assume our covariant star 
products 
to satisfy the following condition.
\begin{de}\label{propB}
Let $(M,\omega,\g,\lambda)$ be a strongly hamiltonian system. Let 
$\star_{\nu}$ be a $\g-$covariant star product on $(M,\omega)$. We 
say 
that $\star_{\nu}$ has the property (B) if there exists an integer 
$N\in\Bbb N$ such that one has
\begin{eqnarray*}
&&(\lambda_{X}\star_{\nu}u)\mod(\nu^{N})=(\lambda_{X}\star_{\nu}u)\mod(\nu^{N+n});\\
&&(u\star_{\nu}\lambda_{X})\mod(\nu^{N})=(u\star_{\nu}\lambda_{X})\mod(\nu^{N+n}),
\end{eqnarray*}
for all $X\in\g$ and $u\in C^{\infty}(M)$ and $n\in\Bbb N$.
\end{de}
In other words the series $\lambda_{X}\star_{\nu}u$ and 
$u\star_{\nu}\lambda_{X}$ stop at order $N$ independently of $u$ in 
$C^{\infty}(M)$.

The data of a covariant star product satisfying the property (B) 
yields representations of $\g$. 
Let
$
E_{\nu}:=C^{\infty}(M)[[\nu,\frac 1{\nu}]]$ be the space of formal 
power series in $\nu$ and $\frac 1{\nu}$ with coefficients in 
$C^{\infty}(M)$.
Assume that the $\g$-covariant star product $\star_{\nu}$ on 
$(M,\omega)$ 
satisfies the property (B).
Then for all $X\in\g$ and $a=\sum_{\ell\in\Bbb 
Z}\nu^{\ell}a_{\ell}\in 
E_{\nu}$, the expression
$$
\lambda_{X}\star_{\nu}a:=\sum_{\ell\in\Bbb 
Z}\nu^{\ell}(\lambda_{X}\star_{\nu}a_{\ell}),
$$
defines an element of $E_{\nu}$.
Indeed, let $c_{k}:C^{\infty}(M)\times C^{\infty}(M)\rightarrow 
C^{\infty}(M)$
be the $k^{th}$ cochain of $\star_{\nu}$, that is
$$
u\star_{\nu}v=:\sum\nu^kc_{k}(u,v),\:\:u,v\in C^{\infty}(M).
$$
Then, $$\lambda_{X}\star_{\nu}a=\sum_{\ell\in\Bbb 
Z}\nu^{\ell}\sum_{k\leq N}\nu^kc_{k}(\lambda_{X},a_{\ell})=
\sum_{m\in\Bbb 
Z}\nu^m\left(\sum_{k+\ell=m}c_{k}(\lambda_{X},a_{\ell})\right).
$$
Therefore, each sum occurring in the parentheses has only a finite 
number of terms since $0\leq k\leq N$.
\begin{de}
Let $(M,\omega,\g,\lambda)$ be a strongly hamiltonian system and 
$\star_{\nu}$ be a $\g$-covariant star product on $(M,\omega)$ 
satisfying
the property (B). One defines the representations 
$\rho^{L}$ 
and $\rho^{R}$ of $\g$ on $E_{\nu}$ by
\begin{eqnarray*}
\rho^{L}(X)a:=&&\frac 1{2\nu}(\lambda_{X}\star_{\nu}a)\:{\rm and}\\
\rho^{R}(X)a:=&&\frac 1{2\nu}(a\star_{\nu}\lambda_{X}).
\end{eqnarray*}
\end{de}
\begin{de}
Let $G$ be a connected Lie group with Lie algebra $\g$. Let 
$(M,\omega,\g,\lambda)$ be a strongly hamiltonian system. Let 
$\star_{\nu}$ be a $\g$-covariant star product on $(M,\omega)$ 
satisfying
the property (B). The associated star representation of $G$ (if it 
exists) is the representation $\pi^{L}$ of $G$ on $E_{\nu}$ such that
$$
d\pi^{L}=\rho^{L}.
$$
\end{de}
          \section{Holonomy reducible symplectic symmetric spaces}
          Let $G$ be a connected simple Lie group. Let us denote by 
$\Lg$ its Lie
          algebra. Let
          $\OO$ be an adjoint orbit of $G$ in $\Lg$. Choose a base 
point $o$ in
          $\OO$ and
          denote by $\Lh$ the Lie algebra of its stabilizer in $\Lg$.
	  Since $\Lg$ is simple the {\it Killing form} 
$\beta$ establishes an equivariant linear 
          isomorphism between $\Lg$ and its dual $\Lg^\star$. 
Therefore, every adjoint orbit 
          can be identified with a coadjoint one. We will denote by 
$\omega^\OO$ the KKS symplectic 
          structure on $\OO$ (cf. Example 1). 

          \begin{dfn} An adjoint orbit $\OO$ in $\Lg$ is {\it 
symmetric} if there
          exists an
          involutive automorphism $\sigma$ of $\Lg$ such that $\Lh=\{ 
X\in \Lg |
          \sigma(X)=X
          \}$. In this case, we denote by $\Lg=\Lh \oplus \Lq$ the 
decomposition of
          $\Lg$
          induced by $\sigma  \, (\sigma|_{\Lq} = -id_{\Lq})$. 
\end{dfn} One has
          $[\Lh,\Lq]=\Lq$ and $ [\Lq,\Lq]=\Lh$.\\
          In a symmetric situation, the KKS form induces on the 
vector
          space
          $\Lq$ a bilinear symplectic form which we denote by 
$\Omega$; this form is
          invariant
          under the action of $\Lh$. The triple $(\Lg, \sigma, 
\Omega)$ is called a
          {\it simple
          symplectic symmetric triple} \cite{[BCG]}. 
	  Symmetric orbits have been studied in \cite{[Bielsss]}. In 
	  particular, on has.
          \begin{prop}\label{p1} Let $\OO$ be a symmetric adjoint 
orbit of a
          simple Lie group $G$ and
          $(\Lg, \sigma, \Omega)$ be its associated symplectic symmetric triple. 
 The following assumptions are
	  equivalent.
	  \begin{enumerate} \item[(i)] the center $\z(\h)$ of $\h$
	      contains a non-compact element.
 \item[(ii)] The subspace $\Lq$ splits 
into a direct
          sum
          $\Lq=\Ll \oplus \Ll'$ of isomorphic $\Lh$-modules. One has 
$[\Ll,\Ll]=0\, ,
          [\Ll',\Ll']=0$ and both $\Ll$ and $\Ll'$ are 
$\beta$-isotropic and
          $\Omega$-Lagrangian subspaces of $\Lq$. 
          \end{enumerate} \end{prop}

           Such a symmetric orbit is called {\it holonomy 
reducible}
          if $\Lh$
          acts reducibly on $\Lq$. 
	  \begin{prop}\label{p2} 
Let $\OO$ be a holonomy reducible symmetric orbit
          in $\Lg$. We define
          the map $ \phi:\Lq=\Ll \oplus \Ll' \to \OO $ by $$
          \phi(l,l'):=Ad(exp(l).exp(l')).o $$ Then $\phi$ is a 
Darboux chart on
          $(\OO,\omega^\OO)$. Precisely, one has
          $\phi^\star\omega^\OO=\Omega$.
          \end{prop}
          We transport the infinitesimal action of $\Lg$ on 
$\phi(\Lq)\subset\OO$ to 
          an infinitesimal action of $\Lg$ on $\Lq$ via $\phi$. One 
then gets a homomorphism 
          of Lie algebra $$\Lg\to\chi(\Lq):X\to X^\star.$$ Setting 
          $$
          \lambda_A(x):=\beta(\phi(x),A) \qquad (x \in \Lq, A\in \Lg),
          $$
          one obtains the strongly hamiltonian system (cf.Definition
	  \ref{strham}) 
          $(\Lq,\Omega,\Lg,\lambda)$. The main property of the 
Darboux chart $\phi$ is
          \begin{prop}
          The Moyal star product on  the symplectic vector space 
$(\Lq,\Omega)$ is 
          $\Lg$-covariant. Moreover the Moyal star product in this 
case satisfies the 
          property (B).
          \end{prop}
          \begin{dfn} Let ${\cal S}_2'(\Lq)$ be the space of 
distributions on $\Lq$
          which are
          tempered in the $\Ll'$-variables w.r.t. the Lebesgue 
measure $dl'$ on
          $\Ll'$. On
          ${\cal S}_2'(\Lq)$, we consider the partial Fourier 
transform ${\cal F}:
          {\cal
          S}_2'(\Lq) \to {\cal S}_2'(\tilde{\Lq})$ which reads 
formally as $$
          ({\cal
          F}u)(l,\eta):=\int_{\Ll'}e^{-i\Omega(\eta, 
l')}u(l,l')dl'\qquad
          i:=\sqrt{-1} $$
          Here $\tilde{\Lq}:=\Ll \oplus \Ll$ i.e. we identify the 
dual space
          ${\Ll'}^\star$
          with $\Ll$ by use of $\Omega$.\\ We will also adopt the 
notation
          $\hat{u}:={\cal F}u$ 
          \end{dfn}
          We define the $\BR$-isomorphism $$
          \tilde{\Lq}\to
          \Ll^{\BC}: (l,\eta) \to z=l +i\nu\eta $$ where the 
parameter $\nu$ is
          now
          considered as being {\it real}.\\
          In
          $\Lq=\Ll\oplus\Ll'$, we
          choose an $\Omega$-symplectic basis $\{ L_a,L_a';\: 1\leq 
a\leq n \}$ where
          $L_a\in\Ll $ and $ L_a'\in \Ll'$.
          On $\Ll$, we define the coordinate system 
$x=\Omega(x,L_a')L_a=:x^aL_a$.\\
          If $h$ is
          an element of $\Lh$, we define its trace as 
	  $$
          \spur(h):=\Omega([h,L_a],L_a').$$
          In this setting, one
          has the
          holomorphic constant vector field on $\Ll^{\BC}$~: $$ 
\partial_{z^a} :=
          \frac{1}{2\nu} (\nu(L_a)_l -i (L_a)_\eta) \quad (1\leq 
a\leq n).$$
          \begin{dfn}\label{ache} Considering $\Ll^{\BC} \subset \Lg^{\BC}$, we 
define, for all
          $A\in
          \Lg^{\BC}$, the polynomials on $\Ll^{\BC}$~: $$ 
\begin{array}{c}
          h_A^{\BC}(z):=
          A_{\Lh^{\BC}} + [A_{{\Ll'}^{\BC}},z] \, \in \Lh^{\BC}\\ 
l_A^{\BC}(z):=
          A_{\Ll^{\BC}}
          + [A_{\Lh^{\BC}},z] +\frac{1}{2}[z,[z,A_{{\Ll'}^{\BC}}]]\, 
\in \Ll^{\BC}
          \end{array}
          $$ where $z\in\Ll^{\BC}$ and
          $A=A_{\Lh^{\BC}}+A_{\Ll^{\BC}}+A_{{\Ll'}^{\BC}}$
          according to the decomposition $\Lg^{\BC}=\Lh^{\BC}\oplus 
\Ll^{\BC} \oplus
          {\Ll'}^{\BC}.$ \end{dfn}

          \begin{dfn}\label{tau} For all $A\in \Lg^{\BC}$, we define the 
holomorphic vector
          field ${\cal
          Z}_A^{(\nu)} \in \Gamma(T^{1,0}(\Ll^{\BC}))$ by $$ 
({\cal
          Z}_A^{(\nu)})_z.f
          :=(l_A^{\BC}(z))^a(\partial_{z^a}.f)(z) $$ where $z\in
          \Ll^{\BC}, f
          \in C^\infty(\Ll^{\BC}, \BC)$ and where, for all 
$w=w_1+i w_2 \in
          \Ll^{\BC}=\Ll
          \oplus i\Ll$, we set $w^a:=w^a_1 +i w^a_2 \qquad (1 
\leq a \leq
          n)$.\\ In the
          same way, extending the Killing form $\beta$ and the trace 
$spur$
          $\BC$-linearly to
          $\Lh^{\BC}$, we define the complex polynomial of degree 1~: 
$$
          \tau^{(\nu)}_A:=\frac{1}{2\nu}(\beta(h_A^{\BC},o)+\nu \, 
spur(h_A^{\BC})).$$
          \end{dfn}\begin{prop} For all $A \in \Lg$, one has $$ 
\frac{1}{2\nu}{\cal
          F}(\lambda_A\star_\nu u)= \tau_A^{(\nu)}.\hat{u} + {\cal
          Z}_A^{(\nu)}.\hat{u} $$
          where $u$ is chosen in such a way that the LHS and the RHS 
make sense (e.g. $u
          \in {\cal S}(\Lq)$ the Schwartz space on $\Lq$). 
          \end{prop}
          In particular, the formula
          \begin{equation}\label{rho}
          \hat{\rho}(A)f:=\tau_A^{(\nu)}f + {\cal
          Z}_A^{(\nu)}.f\qquad(A\in\Lg,f\in 
C^\infty(\l^{\Bbb C})[[\nu,\frac{1}{\nu}]])
          \end{equation}
          defines a (holomorphic) representation of $\Lg$
          on $E_\nu:=C^\infty(\l^{\Bbb C})[[\nu,\frac{1}{\nu}]])$.

\section{Main Theorem}
\subsection{\normalsize Euclidean Jordan algebras and Tube Domains}
An algebra $V$ over $\Bbb R$ or $\Bbb C$ is said to be a Jordan
algebra if for all elements $x$ and $y$ in $V$, one has
\begin{eqnarray*}
x\cdot y&=&y\cdot x,\\ x\cdot(x^{2}\cdot y)&=&x^{2}\cdot(x\cdot y).
\end{eqnarray*}
For an element $x\in V$ let $L(x)$ be the linear map of $V$ defined by
$L(x)y:=x\cdot y$. We denote by $\tau(x,y)$ the symmetric bilinear 
form on $V$ defined by $\tau(x,y)=\Tr L(x\cdot y)$.

A Jordan algebra $V$ is semi-simple if the form $\tau$ is non degenerate 
on $V$. A semi-simple Jordan algebra is unital, we denote by $e$ its 
identity 
element.

One defines on $V$ the triple product
$$
\{x,y,z\}:=(x\cdot y)\cdot z+x\cdot(y\cdot z)-y\cdot (x\cdot z).
$$
We denote by $x\carre y$ the endomorphism of $V$ defined by
$
x\carre y(z):=\{x,y,z\}.
$
Remark that $x \carre y=L(x\cdot y)+[L(x),L(y)]$ See \cite{[Satake]} for more details.
\begin{de}
A Jordan algebra $V_{o}$ over $\Bbb R$ is said to be Euclidean if the bilinear
form $\tau(x,y)$ is positive definite on $V_{o}$.
\end{de}
Let $V_{o}$ be an Euclidean Jordan algebra (EJA), then the set 
$$
C:=\{x^{2}\:\vert\:x\:{\rm invertible}\:{\rm in}\: V_{o}\}
$$
is an open, convex, self-dual cone in $V_{o}$. Those properties of $C$ 
actually characterize $V_{o}$ as an EJA.\newline
Let $V$ be the complexification of $V_{o}$. 
Consider the tube $T_{C}=V_{o}+iC\subset V$
and the Lie group $Aut (T_{C})$ of holomorphic automorphisms of $T_{C}$. 
We denote by $G=(Aut(T_{C}))_{o}$ its identity connected component.

It can be shown that every element $X$ of the Lie algebra $\g$ of the group $Aut(T_{C})$
is a holomorphic vector field on the tube $T_{C}$ of the form (see \cite{[FK]} p.209)
$$
X(z)=u+Tz+P(z)v,
$$
where $u,v\in V_{o}$, $T$ is a linear map of the form $T=a\carre b$ 
with $a,b\in V_{o}$ and $P(z)=2L(z)^2-L(z^{2})$.

In other words, the Lie algebra $\g$
 is a symmetric Lie algebra  
admitting the following graduation
$$
\g=\g_{-1}\oplus\g_{0}\oplus\g_{1},
$$
where $\g_{-1}:=V_{o}$ is the set
of constant polynomial vector fields on $T_{C}$ acting by translations,
$\g_0:=V_{o}\carre V_{o}$ is the subset of 
$\gl(V_{o})$ preserving the cone $C$ and $\g_1:=\{
P(z)v=\{z,v,z\}\:\vert \:v\in V_{o}\}$ is a
subset of homogeneous polynomial maps 
$V\mapsto V$ of degree 2. Remark that $\g_{1}\simeq 
\Ad(j)V_{o}$ where $j\in G$ denotes the Jordan inverse $j(z)=-z^{-1}$.

One writes $(u,T,v)$ for $X\in\g$.
The following result is classical.
\begin{pro}
For $X=(u,T,v)$ and 
$X'=(u',T',v')$ in $\g$
one has
$$
[X,X']=(Tu'-T'u,\:2u'\carre v+[T,T']-2u\carre 
v',\:T'^{\sharp}v-T^{\sharp}v'),
$$
where $T^{\sharp}$ denotes the adjoint endomorphism with respect to 
$\tau$.
\newline
The map $\theta:(u,T,v)\mapsto(v,-T^{\sharp},u)$ is an involutive 
automorphism of $\g$, such that
$\theta(\g_i)=\g_{-i},\:i\in\{-1,0,1\}$. The Lie algebra $\g$ is semi-simple.
\newline
The Killing form on $\g$ is given by
$$
\beta(X,X')=\beta_o(T,T')+2\tr(TT')-4\tau(u,v')-4\tau(v,u'),
$$
where $\beta_o$ denotes the standard Killing form on $\gl(V_{o})$.
\end{pro}
Note that using this identification we have, for $x,y,z\in 
V_{o}=\g_{-1}$
\begin{eqnarray}\label{ident}
x\carre y&=&-\frac 12[x,\theta y]\quad{\rm and}\\ \label{ident2}
\{x,y,z\}&=&-\frac 12\left[[x,\theta y],z\right].
\end{eqnarray}
%\newpage
 The EJA's and their corresponding symmetric Lie algebras 
 are given by the 
following table.
\begin{table}[ht]
 \begin{tabular}{|c|c|c|c|}
\hline
${\g}$ & $\k$ & $V$& $V_{o}$\\
\hline
$\s\u(n,n)$ & $\s\u(n)\oplus i\Bbb R$ & $M(n,\Bbb C)$ &
  $Herm(n,\Bbb C)$\\
$\s\p(n,\Bbb R)$& $\s\u(n)\oplus i\Bbb R$ & $Sym(n,\Bbb C)$&
  $Sym(n,\Bbb R)$\\
$\mathfrak{so}^{*}(4n)$ & $\s\u(2n)\oplus i\Bbb R$ & $Skew(2n,\Bbb 
C)$&
  $Herm(n,\Bbb H)$\\
$\mathfrak{so}(n,2)$ & $\mathfrak{so}(n)\oplus i\Bbb R$ & $\Bbb C^n$& 
$\Bbb R^n$\\
$\e_{7(-25)}$ & $\e_6\oplus i\Bbb R$ & $Herm(3,\Bbb O)\otimes\Bbb C$&
  $Herm(3,\Bbb O)$\\
\hline
\end{tabular}\\
\end{table}

Adopting to our setting results of \cite{[Faraut Olafsson]}, one gets
\begin{pro}
(i) The group $G=Aut(T_{C})_{o}$ admits a symmetric holonomy 
reducible coadjoint orbit.
\newline
(ii) Let $K$ be the maximal compact subgroup of $G$. Then the 
symmetric space $G/K\simeq T_{C}$ is an Hermitian symmetric space of 
tube type.
\end{pro}
\subsection{\normalsize Holomorphic discrete Series}
Let $G=(\Aut(T_{C}))_{o}$ be the Hermitian Lie group associated with 
an Euclidean Jordan algebra $V_{o}$. We denote by $n$ the dimension 
of $V_{o}$ and by $r$ its rank. On $V_{o}$ one defines ( in a 
canonical way) two homogeneous polynomials $\Delta(x)$ of degree $r$
and $\tr(x)$ of degree 1 which coincide with usual determinant and 
trace in case of matrix Jordan algebras. (see \cite{[FK]} for more details).

For a real parameter $m$ consider the space $H^{2}_{m}(T_{C})$ 
of holomorphic functions
$f\in{\cal O}(T_{C})$ such that
$$
\Vert f\Vert^{2}_{m}=\int_{T_C}\vert f(z)\vert^{2}\Delta^{m-2\frac 
nr}(y)dxdy<\infty,
$$
where $z=x+iy\in T_{C}$. Note that the measure $\Delta^{-2\frac 
nr}(y)dxdy$ on $T_{C}$ is invariant under the action of
the group $G$. For $m$ in the 
Wallach set $\cal W$ (see \cite{[FK]} p.264) these spaces are non empty Hilbert 
spaces with reproducing kernels.

The action of $G$ on $H^{2}_{m}(T_{C})$ $(m\in\cal W)$ given by
\begin{equation}
\pi_{m}(g)f(z)=\Det^{m}(D_{g^{-1}}(z))f(g^{-1}.z)
\end{equation}
is called a {\it holomorphic discrete series representation}.

 In the above formula
by $D_{g^{-1}}(z)$ denotes the derivate map of the conformal transformation 
$z\to g^{-1}.z$ of the tube.

The derivate representation 
\begin{equation}\label{sdh}
    d\pi_{m}(X)\phi(z)=-m\frac rn\Tr DX(z)\cdot\phi(z)-D\phi(z)(X(z)).
\end{equation}
 is then given by the following
formul\ae (\cite{[Pev]}).
\begin{eqnarray*}
    {\rm for}\:X(z)=(u,0,0)&&\:d\pi_{m}(X)\phi(z)=-D_{u}\phi(z)\\
  {\rm  
for}\:X(z)=(0,0,v)&&\:d\pi_{m}(X)\phi(z)=-2m\tau(z,v)\phi(z)-D_{p_{v}(z)}
\phi(z)\\
  {\rm for}\:X(z)=(0,T,0)&&\:d\pi_{m}(X)\phi(z)=-m\frac rn\Tr T
  \phi(z)-D_{T(z)}\phi(z),
  \end{eqnarray*}
where $D_{A}\phi(z)$ denotes the action of the tangent vector $A$ on 
the function $\phi$ at the point $z\in V$.
  \begin{thm}\label{resultat}
The star representation $\rho^{L}_{\nu}$ associated to $\star_{\nu}$ is 
equivalent
to the deri\-vate holomorphic discrete series representation $d\pi_m$ 
for $$m=\frac{\beta(o,o)+n\nu c}{4\nu rc}$$
where $c$ is the eigenvalue of the adjoint action of the base point 
$o$.
\end{thm}
{\bf Proof}. 
We use the notations of Section 2.
Using the identification formul\ae $\:$(\ref{ident}) and 
(\ref{ident2}), one gets 
that
$$
{\cal Z}_{X}^{(\nu)}= X(z).
$$
Let us now discuss the expression for 
$
\tau^{(\nu)}_X=\frac 1{2\nu}(\beta(h^{\Bbb C}_X,o)+\nu\spur(h^{\Bbb 
C}_X)).
$
\vskip 5pt

The Lie algebra $\h$ is reductive and therefore it admits the following 
decomposition $\h=\z(h)\oplus[\h,\h]$, where $\z(\h)$ is
the center of $\h$. We write $H=H_{z}+H_{d}$ according to this 
decomposition. Because the trace function vanishes on the Lie algebra 
commutator we have $\spur(H)=\spur(H_{z})$. But the center of $\h$ is
one dimensional (see \ref{p1}) and any element $H_{z}\in\z(\h)$ can be written as
$H_{z}=h_{z}\cdot o$.
Therefore 
$\tau^{(\nu)}_H=\frac 1{2\nu}[\beta(H_z,o)+\nu\spur(H_z)].$ Furthermore we have
$\beta(H_z,o)=h_z\beta(o,o)$ and 
$$
\spur(H_z)=h_z\cdot\ad(o)\vert_{\l}=h_zc\cdot\id\vert_{\l}=h_znc.
$$
In other words it means that $\tau^{\nu}_X$ is proportional to $\Tr(\ad(h_{X}^{\Bbb
C})\vert_{\l})$.

Observe now that for the polynomial $h^{\Bbb C}_X(z)$ introduced in Definition 
\ref{ache}
one has
$
h^{\Bbb C}_X(z)=-2DX(z).
$
So, finally we have
$$
\tau^{(\nu)}_X=-\frac{\beta(o,o)+n\nu c}{4n\nu c}\Tr DX.
$$
The identification of corresponding terms in formul\ae$\:$ (\ref{rho})
and (\ref{sdh}) completes the proof.\qed
{\bf Remark}. The
vector
field part ${\cal Z}_X^{(\nu)}$ of the star representation
$\hat\rho(X)$ $X\in\g$ ( see (\ref{rho})) is in general
singular on the
entire $V=\l^{\Bbb C}$. By this we mean that, denoting by $\phi_t^X$ the local
flow of $X$, the set of ``bad" Cauchy datas
$S_{X}:=\{z\in V\, |\, \phi_t^X(z)\mbox{ is not defined for all values of
}t\,\}$ is in general not empty. However, in the case where the group $G$
is the automorphism group of a tube domain, Theorem \ref{resultat} shows that
the complementary set $U$ of $\cup_{X\in\g}S_X=V_o$ in $V$ is not empty. The infinitesimal action
of $\g$ on (a connected component of) $U$ exponentiates to $G$ as its
action on the holomorphic tubular realization of $G/K$. A deep geometric
study and in particular relations between star representations and Riccati type ordinary differential equations
(see \cite{[Zel]}) has not been investigated here. However, one can at least say that, given a holonomy
reducible symmetric
orbit ${\cal O}$ of $G$, one then canonically gets a one parameter family
of representations $\{\rho_\nu\}_{\nu\in\R}$ of $\g$ deforming the
infinitesimal action ($\rho_0$) of $\g$ on an open subset $\phi(\p)$ of ${\cal
O}$. This parameter family leads to an interpolation between
the infinitesimal action of $\g$ on $G/H$ and its holomorphic action on
$G/K$ ($\rho_\nu$ for the value $-\frac{1}{nc}\beta(o,o)$ of $\nu$).

\end{document}